\documentstyle[12pt]{article}
\newtheorem{Th}{\hspace*{\parindent} Theorem}
\newtheorem{Prop}[Th]{\hspace*{\parindent} Proposition}

\newtheorem{Rem}[Th]{\hspace*{\parindent} {\it Remark}}
\newtheorem{Cor}[Th]{\hspace*{\parindent} Corollary}
\newtheorem{Lemma}[Th]{\hspace*{\parindent} Lemma}
\newcommand{\Proof}{{\it Proof. }}

\newcommand{\be}{\begin{eqnarray*}}
\newcommand{\ee}{\end{eqnarray*}}
\newcommand{\pkEF}{{\cal P}(^k\!E,F)}

\newcommand{\pkE}{{\cal P}(^k\!E)}

\newcommand{\pwbkEF}{{\cal P}_{wb}(^k\!E,F)}
\newcommand{\pwbkE}{{\cal P}_{wb}(^k\!E)}
\newcommand{\pwcokEF}{{\cal P}_{wco}(^k\!E,F)}
\newcommand{\pcckEF}{{\cal P}_{cc}(^k\!E,F)}
\newcommand{\pcck}{{\cal P}_{cc}(^k\!}

\newcommand{\Co}{{\cal C}o\,}
\newcommand{\WCo}{{\cal WC}o\,}
\newcommand{\CC}{{\cal CC}}

\newcommand{\LEF}{{\cal L}(E,F)}

\newcommand{\pwcok}{{\cal P}_{wco}(^k\!}

\newcommand{\finesp}{\hspace*{\fill} $\Box$\vspace{0.5\baselineskip}}
\newcommand{\fin}{\hspace*{\fill} $\Box$}
\newcommand{\pol}{{\cal P}(}
\newcommand{\ptp}{\widehat{\bigotimes}} 
\newcommand{\sptp}{\widehat{\Delta}} 
\newcommand{\Gr}{{\cal G}r}

\newcommand{\lra}{\longrightarrow}
\newcommand{\Ra}{\Rightarrow}

\newcommand{\ra}{\rightarrow}

\newcommand{\bmi}{\noindent \begin{minipage}[t]{136mm}}

\newcommand{\emi}{\end{minipage}}

\newcommand{\ssc}{\subsection*}

\begin{document}

\title{Polynomial Grothendieck properties
      }

\author{Manuel Gonz\'alez\thanks{Supported in part by DGICYT Grant PB
      91--0307 (Spain)} \\
      Departamento de Matem\'aticas \\
      Facultad de Ciencias\\
      Universidad de Cantabria \\ 39071 Santander (Spain)
       \and
      Joaqu\'\i n  M. Guti\'errez\thanks{Supported in part by DGICYT
      Grants PB 90--0044 and PB 91-0307 (Spain)}\\
      Departamento de Matem\'atica Aplicada\\
      ETS de Ingenieros Industriales \\
      Universidad Polit\'ecnica de Madrid\\
      C. Jos\'e Guti\'errez Abascal 2 \\
      28006 Madrid (Spain)}

\date{\mbox{}}
\maketitle

\begin{abstract}
A Banach space $E$ has the Grothendieck property if every (linear bounded)
operator from $E$ into $c_0$ is weakly compact.
It is proved that, for an integer $k>1$,
every $k$-homogeneous polynomial from $E$
into $c_0$ is weakly compact if and only if the space
${\cal P}(^kE)$ of scalar
valued polynomials on $E$ is reflexive.
This is equivalent to the symmetric $k$-fold projective tensor product
of $E$ (i.e., the predual of ${\cal P}(^kE)$) having the Grothendieck property.
The Grothendieck property of the projective tensor product
$E\widehat{\bigotimes}F$
is
also characterized. Moreover,
the Grothendieck property of $E$ is described in terms of sequences of
polynomials.
Finally, it is shown that if every operator from $E$ into $c_0$
is completely continuous, then so is every polynomial between these
spaces.
\end{abstract}

Throughout, $E$, $F$ will be Banach spaces, and $E^*$ the dual
of $E$.
We denote by $\LEF$ the space of all (linear bounded) operators from $E$
to $F$, and by $\Co (E,F)$ ($\WCo (E,F)$) the subspace of all
(weakly) compact operators.
We say that $T \in \LEF$ is {\it completely continuous\/} if it takes weakly
convergent sequences into norm convergent sequences, and we write
$T\in\CC (E,F)$.

For an integer $k$, we shall consider the following classes of polynomials:

(a) $\pkEF$ is the space of all $k$-homogeneous (continuous) polynomials from
      $E$ to $F$;

(b)
$\pcckEF$, the subspace of {\it completely continuous
polynomials,} i.e., the polynomials taking weakly convergent sequences into norm
convergent ones, equivalently, taking weak Cauchy sequences
into convergent ones \cite[Theorem~2.3]{AHV};

(c) $\pwcokEF$, the subspace of weakly compact polynomials.

(d) $\pwbkEF$, the polynomials whose restrictions to bounded subsets of $E$
      are weakly continuous; these are compact polynomials. It is well kown
      that $\pwbkEF \subseteq \pcckEF$, and the equality occurs if and only if
      $E$ contains no copy of $\ell_1$ (see e.g.\ \cite{Gu}).

      The space of $k$-linear (continuous) mappings from $E^k$ to $F$ is denoted
by ${\cal L}(^k\!E,F)$. To each $P\in\pkEF$ we can associate a unique symmetric
$A\in {\cal L}(^k\!E,F)$ such that $P(x)=A(x,\ldots ,x)$ for all $x\in E$.
Whenever $F$ is omitted, it is understood to be the scalar field.
For the general theory of polynomials between Banach spaces, we refer to
\cite{Mu}.

      The projective tensor product of $E$ and $F$ is referred to as
$E\ptp F$. The closed linear span of the set $\{ x\otimes
\stackrel{(k)}{\ldots} \otimes x: x\in E\}$ in $\ptp ^k E := E\ptp
\stackrel{(k)}{\ldots}\ptp E$ is denoted by
$\sptp ^kE$. Its dual is isomorphic to $\pkE$. The spaces $\pkEF$ and
${\cal L}(\sptp^kE,F)$ are linearly isometric, and the image of
$\pwcokEF$ under this isometry is $\WCo (\sptp^kE,F)$ \cite{RyTh}.

We say that $E$ has the {\it Grothendieck property,} and write $E\in \Gr$,
if every sequence in $E^*$ converging to zero in the weak-star $(w^*)$
topology, is also weakly null. Equivalently, if every operator $E\ra c_0$
is weakly compact.

In this paper, we investigate conditions on $E$ so that every $k$-homogeneous
polynomial $E\ra c_0$ be weakly compact, equivalently $\sptp^kE \in \Gr$,
proving that this is the case if and only if $\pkE$, and hence $\sptp^kE$, is
reflexive.
These ``polynomially reflexive'' Banach spaces have been investigated by
various authors \cite{AAD,Fa,GGUC,GJ}.

      We also show that the situation is different if we consider the
Grothendieck property of $E\ptp F$, proving that for $E\in \Gr$ and
$F^*$ reflexive with the bounded compact approximation property,
$E\ptp F \in \Gr$ if and only if ${\cal L}(E,F^*) = \Co (E,F^*)$.
In particular,
$\ell_\infty \ptp \ell_p \in \Gr$ for $2<p<\infty$, and
$\ell_\infty \ptp T^* \in\Gr$, where $T^*$ is the original Tsirelson space.

      Every $P\in \pkE$ has a standard extension $\tilde{P}\in \pol ^k\!
E^{**})$
(see \cite{DG}).
Thus, given $P\in \pol ^k\!E, c_0)$, with
$Px= (P_nx)_n$, we can define $\tilde{P}\in \pol ^k\!E^{**},\ell_\infty)$
by $\tilde{P}z:= (\tilde{P}_nz)_n$.
It is shown that $E\in \Gr$ if and only if, for every sequence $(P_n) \subset
\pkE$ such that $P_nx\ra 0$ for all $x\in E$, we have $\tilde{P}_nz\ra 0$
for every $z\in E^{**}$.
Then $E\in\Gr$ if and only if for every $P\in \pol ^k\! E, c_0)$,
we have $\tilde{P}(E^{**}) \subseteq c_0$.

      Several authors \cite{Pel,GJ,GGwh} have studied conditions on $E$, $F$
so that $\pkEF = \pcckEF$. Here we investigate the equality
$\pol ^k\! E,c_0) = \pcck E, c_0)$, proving that it is equivalent to
${\cal L}(E,c_0) = \CC (E, c_0)$.
Grothendieck spaces with the Dunford-Pettis property, and
Schur spaces satisfy this property.
In \cite{L}, examples are given of Grothendieck spaces with the Dunford-Pettis
property, and semigroups of operators on these spaces are studied.

\ssc{\S 1}

      We first characterize the spaces $E$ such that
$\pol ^k\!E,c_0) = \pwcok E,c_0)$.
Some previous results are needed.

\begin{Prop}
\label{GrnoR}
\hspace*{1em}\\
\hspace*{\parindent}
Let $E$ be a nonreflexive space with the Grothendieck property. Then $E$
contains a copy of $\ell_1$.
\end{Prop}

\Proof

Since $E\in\Gr$, $E$ has no quotient isomorphic to $c_0$.
Assume $E$ contains no copy of $\ell_1$. Then, $E^*$ contains no copy
of $\ell_1$ \cite[Corollary~2.3]{GO}.
Given a bounded sequence $(\phi_n) \subset E^*$, we can
find a weak Cauchy subsequence $\left( \phi_{n_k}\right)$. The sequence
$\left( \phi_{n_k}\right)$ is $w^*$-convergent, hence weakly convergent,
and we conclude that $E$ is reflexive.
\fin

\begin{Prop}
\label{allwco}
\hspace*{1em}\\
\hspace*{\parindent}
If $\pol ^k\! E,c_0) = \pwcok E,c_0)$ for some $k > 1$, then $E$ is reflexive.
\end{Prop}

\Proof

Suppose there is a nonweakly compact $T\in {\cal L}(E,c_0)$.
Then we can find a bounded sequence $(x_n) \subset E$ such that $(Tx_n)$
converges in the topology
$\sigma (\ell_\infty ,\ell_1)$ to some $(a_n)\in\ell_\infty$
with $\lim a_n = a \neq
0$. Define $P(x): = (Tx)^k$, i.e., take the $k\,$th power coordinatewise.
We have $P\in \pwcok E,c_0)$. However, $(Px_n)$ converges
in the topology $\sigma (\ell_\infty ,\ell_1)$ to the sequence $(a_n^k)_n \in
\ell_\infty$,
with $\lim_n a_n^k = a^k \neq 0$, a contradiction.
Hence, $E\in\Gr$.

      Suppose $E$ is nonreflexive. By Proposition~\ref{GrnoR}, $E$
contains a copy of $\ell_1$. Then, there is a quotient map $q:E\ra \ell_2$
(see e.g.\ \cite[Lemma~12]{GGUC}).
Let $P:\ell_2 \ra \ell_1$ be the polynomial given by
$P((x_n)_n) = (x_n^k)_n$, and let $q':\ell_1\ra c_0$ be a quotient map.
Then the product $q'Pq \in \pol ^k\!E,c_0)$ is not weakly compact, a
contradiction.
\fin

\begin{Lemma}
\label{compl}
\hspace*{1em}\\
\hspace*{\parindent}
If $M$ is a complemented subspace of $E$, then $\sptp^kM$ is a complemented
subspace of $\sptp^kE$.
\end{Lemma}

\Proof

If $S\in {\cal L}(E,E)$ is a projection with $S(E)=M$, consider the linear
mapping defined by
$$
x\otimes \stackrel{(k)}{\ldots} \otimes x \longmapsto
Sx \otimes \stackrel{(k)}{\ldots}\otimes Sx \hspace{2em} (x\in E)\, .
$$
Easily, this mapping extends to an operator in ${\cal L}(
\sptp^kE ,\sptp^kE)$, which is the required projection.
\finesp

      To simplify notation, we write $x^{(k)}:=x\otimes\stackrel{(k)}{\ldots}
\otimes x$, for $x\in E$.
Let $P_k: \otimes^k E\ra \otimes ^kE$ be the projection defined by
$$
P_k(x_1\otimes\cdots\otimes x_k) = \frac{1}{k!2^k}\sum_{\epsilon_j=\pm 1}
\epsilon_1 \cdots\epsilon_k \left( \epsilon_1x_1 + \cdots +\epsilon_kx_k\right)
^{(k)}\, .
$$
This mapping extends to an operator on $\ptp ^kE$, which is a projection of
$\ptp^kE$ onto $\sptp ^kE$.
The following Lemma is contained in \cite{Go}. We include the proof
for completeness.

\begin{Lemma}
\label{symtpr}
\hspace*{1em}\\
\hspace*{\parindent}
Given $u\in\sptp ^kE$, there exists a sequence $(x_i)\subset E$ such that
$\sum_{i=1}^\infty \| x_i\| ^k<\infty$ and $u=\sum_{i=1}^\infty\epsilon_i
x_i^{(k)}$, with $\epsilon_i\in \{ \pm 1\}$.
\end{Lemma}

\Proof

By the definition of the projective norm, for any $\delta>0$, we can find
sequences
$\left( x_n^1\right) ,\ldots ,\left( x_n^k\right) \subset E$, such that
$$
\sum_{n=1}^\infty \left\| x_n^1\right\|\cdot \ldots \cdot\left\| x_n^k\right\|
< \| u\| +\delta\, , \mbox{ and } u = \sum_{n=1}^\infty x_n^1\otimes \cdots
\otimes x_n^k\, .
$$
Then,
$$
u = P_ku = \sum_{n=1}^\infty P_k\left( x_n^1\otimes \cdots \otimes
x_n^k\right)\, .
$$
We can assume that, for each $n$, $\left\| x_n^1\right\| = \cdots = \left\|
x_n^k\right\|$. Since
$$
P_k\left( x_n^1\otimes \cdots\otimes x_n^k\right) = \frac{1}{k!2^k}
\sum_{\epsilon_j=\pm 1} \epsilon_1\cdots\epsilon_k
\left( \epsilon_1 x_n^1 +\cdots +\epsilon_kx_n^k\right) ^{(k)}\, ,
$$
denoting by $(x_i)$ the sequence
$$
\left\{ \frac{1}{2\sqrt[k]{k!}}
\left( \epsilon_1 x_n^1 +\cdots +\epsilon_kx_n^k\right) :
\epsilon_1,\ldots ,\epsilon_k=\pm 1; \;
n=1,2,\ldots \right\} \, ,
$$
we obtain
$$
\sum_{i=1}^\infty \left\| x_i\right\| ^k \leq \frac{k^k}{k!} \left( \| u\| +
\delta\right) \, , \mbox{ and }
u=\sum_{i=1}^\infty \epsilon_ix_i^{(k)} \, ,
$$
completing the proof.
\finesp

      Before stating the main result of this part, recall that for $P\in\pkEF$,
its {\it adjoint\/} is the operator $P^*:F^*\ra\pkE$ given by
$P^*(\psi) = \psi
\circ P$ for every $\psi \in F^*$.
Then $P\in\pwcokEF$ if and only if $P^*$ is weakly compact
\cite[Proposition~2.1]{RyWC}.

\begin{Th}
\label{polGr}
\hspace*{1em}\\
\hspace*{\parindent}
Given $k > 1$, we have that $\pol ^k\! E,c_0) = \pwcok E,c_0)$ if and only
if the space $\pkE$ is reflexive.
\end{Th}

\Proof

For the ``only if'' part, if $E$ is separable, then $\sptp^kE$ is a Grothendieck
separable space, hence reflexive. In the general case, suppose $\sptp^kE$ is
not reflexive. By Lemma~\ref{symtpr}, we can find
$$
w_n = \sum_{i=1}^\infty \epsilon_i
x_n^i \otimes \stackrel{(k)}{\ldots} \otimes x_n^i
\; , \hspace{1.5em} \| w_n\| \leq 1
$$
so that $\{ w_n\}$ is not relatively weakly compact in $\sptp^kE$. Let
$N$ be the closed linear span of $\{x_n^i\} _{i,n}$ in $E$. Since $E$ is
reflexive (Proposition~\ref{allwco}),
there is a separable subspace $M$ complemented
in $E$ with $N\subseteq M$.
By Lemma~\ref{compl}, $\sptp^kM \in \Gr$, and is separable, therefore reflexive.
However, $(w_n)\subset \sptp^kM$, a contradiction, and we conclude that
$\sptp ^k E$ and $\pkE$ are reflexive.

      The ``if'' part is clear by the previous comment.
\finesp

      It is well known that the space $\pol ^k\ell_p)$ is reflexive if and
only if $k<p<\infty$. If $E=T^*$, then $\pkE$ is
reflexive for all $k$ \cite{AAD}.

      In the last Theorem, $c_0$ can be replaced by any superspace $F$. However,
for $F$ containing no copy of $c_0$, the result is not true since, for instance,
every polynomial from $E=c_0$ into $F\not\supset c_0$
is weakly continuous on bounded subsets
(see e.g.\ \cite{GGwh}).

\ssc{\S 2}

      In this part, we show that
the situation is different for general projective tensor products.
We refine a result of \cite{Kh}, proving
that for $E\in\Gr$ and $F^*$ reflexive with the bounded compact
approximation property, $E\ptp F\in\Gr$ if and only if ${\cal L}(E,F^*) =
\Co (E,F^*)$.
As a consequence, $\ell_\infty \ptp \ell_p$ for $2<p<\infty$ and $\ell_\infty
\ptp T^*$ have the Grothendieck property.

\begin{Prop}
\label{oneref}
\hspace*{1em}\\
\hspace*{\parindent}
Suppose $E\ptp F \in \Gr$. Then $E, F\in\Gr$ and at least one of them is
reflexive.
\end{Prop}

\Proof

Since $E$ and $F$ are complemented in $E\ptp F$, the first assertion is clear.
Suppose $E$ and $F$ are nonreflexive. Then each of them contains a copy of
$\ell_1$ (Proposition~\ref{GrnoR}). Hence, there are quotient maps (see e.g.\
\cite[Lemma~12]{GGUC})
$$
q_1:E\lra \ell_2 \hspace{1em}\mbox{ and }\hspace{1em} q_2:F\lra \ell_2 \, .
$$
Consider the quotient maps
$$
E\ptp F \stackrel{q_1\ptp \mbox{\small id}}{\lra} \ell_2 \ptp F
\stackrel{\mbox{\small id}\ptp q_2}{\lra} \ell_2 \ptp \ell_2 \, .
$$
It is well known that $\ell_2 \ptp \ell_2 \not\in \Gr$ (separable Gronthendieck
spaces are reflexive). Hence, $E\ptp F \not\in \Gr$, a contradiction.
\finesp

\begin{Rem}
\hspace*{1em}\\
\hspace*{\parindent}
{\rm It follows from Proposition~\ref{oneref} that whenever
$E_1 \ptp \cdots \ptp E_k \in \Gr$ and, for example, $E_1$ is not reflexive,
$E_2 \ptp \cdots \ptp E_k$ is reflexive. In particular, $E_2,\ldots ,E_k$ are
reflexive.}
\end{Rem}

      Before stating the next result, recall that the dual of $E\ptp F$ may
be identified with ${\cal L}(E,F^*)$.

\begin{Prop}
\label{allcoimGr}
\hspace*{1em}\\
\hspace*{\parindent}
Assume $E\in\Gr$ and $F$ is reflexive. If ${\cal L}(E,F^*) = \Co (E,F^*)$,
then $E\ptp F\in\Gr$.
\end{Prop}

\Proof

Let $(A_n)\subset {\cal L}(E,F^*)$ be a $w^*$-null sequence.
Then for every $x\in E$ and $y \in F$,
$$
\langle y, A_n(x)\rangle = \langle x\otimes y,A_n\rangle \lra 0 \, .
$$
Applying Kalton's test for the weak convergence of sequences in spaces of
compact operators (see Theorem~1 or 3 in \cite{Ka}),
we have that $(A_n)$ is weakly null.
\finesp

      We say that $E$ has the {\it bounded compact approximation property\/}
(BCAP) \cite{AT} if there exists $\lambda \geq 1$ so that for each compact
subset $K\subset E$ and for each $\epsilon >0$ there is $S\in \Co (E,E)$
such that
$$
\sup \left\{ \|Sx - x\| : x\in K\right\} \leq \epsilon
\mbox{ , \hspace{1.5em}} \| S - \mbox{id}\| \leq \lambda \, .
$$
Every space with the bounded approximation property has the BCAP.
The converse is not true \cite{W}.

\begin{Prop}
\label{Grimallco}
\hspace*{1em}\\
\hspace*{\parindent}
Suppose $F^*$ is reflexive and has the {\rm BCAP}, and $E\ptp F\in\Gr$.
Then we have ${\cal L}(E,F^*) = \Co (E,F^*)$.
\end{Prop}

\Proof

Suppose first that $F^*$ is separable. Then there is a bounded sequence
$(T_n) \subset \Co (F^*,F^*)$ such that $T_n\psi \ra \psi$ for all
$\psi \in F^*$.
Assume $T\in {\cal L}(E,F^*)$ is not compact. For $x\in E$, $y\in F$
we have
$$
\langle x\otimes y, T_nT\rangle = \langle y,T_n(Tx)\rangle \lra
\langle y,Tx\rangle = \langle x\otimes y, T\rangle \, .
$$
Since $(T_nT)$ is bounded and $\{ x\otimes y: x\in E,\, y \in F\}$ generates
a dense subset of $E\ptp F$, we have that $(T_nT)$ is $w^*$ convergent to
$T$. Since $(T_nT)\subset \Co (E, F^*)$, $(T_nT)$ is not weakly convergent
to $T$, a contradiction.

      For $F^*$ nonseparable, suppose $T$ as above. There is a bounded sequence
$(x_n)\subset E$ such
that $(Tx_n)$ has no Cauchy subsequence. The closed linear span of
$\{ Tx_n\}$ is contained in a separable space $M^*$ complemented in $F^*$.
If $q:F^* \ra M^*$ is the identity on $M^*$, then $qT\in {\cal L}(E,M^*)$ is
noncompact.
By the above, $E\ptp M\not\in \Gr$, a contradiction since $E\ptp M$
is a quotient of $E\ptp F$.
\fin

\begin{Cor}
\hspace*{1em}\\
\hspace*{\parindent}
For $1\leq p\leq\infty$,
the space $\ell_\infty \ptp \ell_p$ has the
Grothendieck property if and only if $2<p<\infty$.
\end{Cor}

\Proof

If $2<p<\infty$, we have $\ell_p^* = \ell_q$ with $1<q<2$, and it is known
that every operator $\ell_\infty \ra \ell_q$ factors through $\ell_2$
\cite[Corollary~4.4]{Pi} and is therefore compact.
The converse is easy.
\fin

\begin{Rem}
\hspace*{1em}\\
\hspace*{\parindent}
{\rm a) We note that   the space
$\ell_\infty \ptp \ell_p\ptp\ell_p$ does not have the Grothendieck
property, for $2\leq p\leq 3$.
Indeed, there is a noncompact operator $T:\ell_p\ra \left( \ell_p \ptp\ell_p
\right) ^*$, for instance, the operator $T$ associated to the polynomial
$Px:=\sum_{i=1}^\infty x_i^3$, for $x=(x_i)\in\ell_p$,
given by $(Tx)(y\otimes z)=\hat{P}(x,y,z)$, where $\hat{P}$ is the
symmetric 3-linear form associated to $P$.
Then $T$ is not completely continuous, so we can find a weakly null
sequence $(x_n)\subset\ell_p$ such that $\{ Tx_n\}$ is not relatively
compact. Passing to a subsequence, we can assume that $(x_n)$ is equivalent
to a block basis and hence equivalent to the basis of $\ell_p$.
Let $q:\ell_\infty\ra\ell_2$ be a quotient, and
$j:\ell_2\ra\ell_p$ the operator taking the $\ell_2$-basis into $(x_n)$.
Then $Tjq:\ell_\infty\ra\left( \ell_p\ptp\ell_p\right) ^*$ is not compact,
and it is enough to apply Proposition~\ref{Grimallco}.

b) It is proved in \cite[Corollary~8]{AAF} that the space $\left( \ptp ^kT^*
\right) \ptp\ell_p$ is reflexive, for $1<p<\infty$.}
\end{Rem}

\begin{Cor}
\hspace*{1em}\\
\hspace*{\parindent}
The space $\ell_\infty \ptp T^*$ has the Grothendieck property.
\end{Cor}

The proof relies on the following Lemma.

\begin{Lemma}
\hspace*{1em}\\
\hspace*{\parindent}
Let $T$ be the dual of $T^*$.
Then ${\cal L}(\ell_\infty,T) = \Co (\ell_\infty,T)$.
\end{Lemma}

\Proof

Assume $S\in {\cal L}(\ell_\infty,T)$ is not compact.
Let $(y_k)\subset\ell_\infty$ be a bounded sequence such that $(Sy_k)$ has
no convergent subsequence.
Choose a weakly convergent subsequence $\left( Sy_{k_n}\right)$, and take
$x_n:= y_{k_{2n}} - y_{k_{2n-1}}$.
Then $(Sx_n)$ is weakly null, and we can assume that it is equivalent to a block
basis in $T$.

      Since $\{ Sx_n\}$ spans a complemented subspace $[Sx_n]$
\cite[Proposition~II.6]{CS}, there is an operator $V:T\ra [Sx_n]$
which is the identity on $[Sx_n]$.
For $1<q<2$, $T$ has lower $q$-estimates \cite[Proposition~V.10]{CS},
so there is an operator $U:[Sx_n]\ra\ell_q$ given by $U(Sx_n)=e_n$, where
$(e_n)$ is the unit vector basis of $\ell_q$.
Then $UVS:\ell_\infty\ra\ell_q$ is not compact, a contradiction
\cite[Corollary~4.4]{Pi}.
\finesp

      In
\cite{Kh}, the following result was obtained
(see {\it Zentralblatt Math.} {\bf 599} \#46017
(1987)):

``Let $E$ be a Banach space with the Grothendieck property, and $F$ a reflexive
space with the metric approximation property.
For $E\ptp F$ to have the Grothendieck property it is necessary and sufficient,
that each operator $E\ra F^*$ be compact.''

      Another related result is the following of \cite{H}:

      ``If $E$ and $F$ are reflexive and both have the approximation property,
then $\LEF$ is reflexive if and only if $\LEF = \Co (E,F)$.''

\ssc{\S 3}

      Next we describe the Grothendieck property in terms of polynomials.
Recall that each $P\in\pkE$ has a Davie-Gamelin extension $\tilde{P} \in
\pol ^k E^{**})$ (see the Introduction).
The authors are indebted to Professor Richard M. Aron, who suggested this study.
Namely, he asked if, given $E\in\Gr$ and a sequence $(P_n)\subset \pkE$ with
$P_nx\ra 0$ for all $x\in E$, it is true that $\tilde{P}_nz\ra 0$ for all
$z\in E^{**}$.
The following Theorem shows that the answer is affirmative.

\begin{Th}
\label{Grseqpol}
\hspace*{1em}\\
\hspace*{\parindent}
The following assertions are equivalent:

{\rm (a)} $E$ has the Grothendieck property;

{\rm (b)} for every integer $k$, given a sequence $(P_n) \subset \pkE$ with
      $P_nx\ra 0$ for all $x\in E$, then $\tilde{P}_nz\ra 0$ for all
      $z\in E^{**}$;

{\rm (c)} the same statement as {\rm (b)} is true for some $k$.
\end{Th}

\Proof

(a) $\Ra$ (b)
By induction on $k$. For $k=1$, the result is nothing but the definition
of the Grothendieck property.
Suppose it holds for
$k-1$, and let $(P_n)\subset \pkE$ be a sequence such that $P_nx\ra 0$ for
all $x\in E$.
Denote by $F_n\in {\cal L}(^k\!E)$ the associated symmetric $k$-linear
form, and by
$G_n\in {\cal L}\left(
E\times E^{**}\times \stackrel{(k-1)}{\ldots}\times E^{**}\right)$
an extension obtained by the Davie-Gamelin method.

Thanks to the polarization formula \cite[Theorem~1.10]{Mu}, we have that
$$
F_n(x_1,x_2,\ldots ,x_k) \lra 0 \hspace{1em}\mbox{ for every }\hspace{1em}
x_1,\ldots ,x_k\in E\, .
$$
Fixing $x_1\in E$, we define $Q_n\in {\cal P}(^{k-1}\!E)$ by
$Q_n(x)= F_n(x_1,x,\ldots ,x)$.
Then $Q_nx\ra 0$ for all $x\in E$.

By the induction hypothesis and polarization,
$$
G_n(x_1,z_2,\ldots ,z_k)\lra 0\, ,\hspace{1em}\mbox{ for }\hspace{1em}
z_2,\ldots ,z_k\in E^{**}\, .
$$
Then, for $z_2,\ldots ,z_k\in E^{**}$ fixed, the sequence $(\phi_n)\subset
E^*$, given by $\phi_n(x) = G_n(x,z_2,\ldots ,z_k)$, is $w^*$-null,
hence weakly null, and so
$$
\tilde{F}_n (z_1,z_2,\ldots ,z_k)\lra 0 \hspace{1em}\mbox{ for every }
\hspace{1em}
z_1,\ldots ,z_k\in E^{**} \, ,
$$
where $\tilde{F}_n\in {\cal L}(^k\!E^{**})$ is the Davie-Gamelin extension of
$F_n$.

(b) $\Ra$ (c) is trivial.

(c) $\Ra$ (a) Given a $w^*$ null sequence $(\phi_n)\subset E^*$, apply (c)
to $P_nx: = \left( \phi_n(x)\right) ^k$.
\finesp

      Given a polynomial $P\in {\cal P}(^k\!E,c_0)$, with $Px=(P_nx)_n$,
we define $\tilde{P} \in {\cal P}(^k\!E^{**},\ell_\infty )$ by
$\tilde{P}z:= (\tilde{P}_nz)_n$. Then we conclude

\begin{Cor}
\hspace*{1em}\\
\hspace*{\parindent}
The space $E$ has the Grothendieck property if and only if for every
$P\in {\cal P}(^kE,c_0)$, we have that $\tilde{P}(E^{**})\subseteq c_0$.
\end{Cor}

      For polynomials whose restrictions to bounded sets are weakly continuous,
we can deduce a result on weak convergence:

\begin{Cor}
\hspace*{1em}\\
\hspace*{\parindent}
The following assertions are equivalent:

{\rm (a)} $E$ has the Grothendieck property;

{\rm (b)} for every integer $k$ and every $F$, if for a sequence $(P_n)\subset
      \pwbkEF$ we have that $\langle P_nx,\psi\rangle\ra 0$ for all $x\in E$
      and $\psi\in F^*$, then $(P_n)$ is weakly null;

{\rm (c)} the same statement as {\rm (b)} is true for some $k$ and some
      $F\neq \{ 0\}$;

{\rm (d)} for some $k$, if for a sequence $(P_n)\subset \pwbkE$ we have that
      $P_nx\ra 0$ for all $x\in E$, then $(P_n)$ is weakly null.
\end{Cor}

\Proof

(a) $\Ra$ (b) It is proved in Theorem 4 of \cite{GGwc} that a sequence
$(P_n)\subset \pwbkEF$ is weakly null if and only if, for every $z\in E^{**}$
and $\psi \in F^*$, we have $\langle \tilde{P}_nz,\psi\rangle\ra 0$.
Therefore, it is enough to apply Theorem~\ref{Grseqpol}(b).

(b) $\Ra$ (c) is trivial.

(c) $\Ra$ (d) Take $0\neq y\in F$ and define $Q_nx: = (P_nx)y$.

(d) $\Ra$ (a) Given a $w^*$-null sequence $(\phi_n)\subset E^*$, apply (d)
to $P_nx: = \left( \phi_n(x)\right) ^k$.
\fin

\ssc{\S 4}

      Several authors \cite{Pel,GJ,GGwh} have studied conditions on $E$, $F$
so that $\pkEF = \pcckEF$. Here we investigate the equality
$\pol ^k\! E,c_0) = \pcck E, c_0)$, proving that it is equivalent to
${\cal L}(E,c_0) = \CC (E, c_0)$.
Therefore, the Grothendieck spaces with the Dunford-Pettis property, and the
Schur spaces satisfy this property.

\begin{Th}
\hspace*{1em}\\
\hspace*{\parindent}
The following assertions are equivalent:

{\rm (a)} ${\cal L}(E,c_0) = \CC (E,c_0)$;

{\rm (b)} ${\cal P}(^k\!E,c_0) = \pcck E,c_0)$ for all integers $k$.

{\rm (c)} ${\cal P}(^k\!E,c_0) = \pcck E,c_0)$ for some integer $k$.
\end{Th}

\Proof

(a) $\Ra$ (b) By induction on $k$. For $k=1$ there is nothing to prove.
Assume the result
is true for $k-1$, and consider $P\in {\cal P}(^k\!E,c_0)$ with
associated $k$-linear mapping $A$.
We only sketch the proof, since it  follows the lines of that
in \cite[Theorem~6]{GGwh}.
In fact, we can prove that every $k$-linear mapping from $E^k$ into $c_0$
takes weak Cauchy sequences into convergent ones.
Let $\left( x^1_n\right) , \ldots ,
\left( x_n^k\right) \subset E$ be weak Cauchy sequences. Suppose first
that one of them, say $\left( x^1_n\right)$, is weakly null.

Define the operator
\be
T:E & \lra & c_0(c_0) \\
y   & \longmapsto & \left( A\left( x^1_n, \ldots x_n^{k-1},y\right)\right) _n
\, .
\ee
Using the induction hypothesis,
it is not difficult to  see that $T$ is well-defined.
Since $c_0(c_0)$ is isomorphic to $c_0$,
$T$ is completely continuous.
>From this, we have
$$
\lim_m \left\| A\left( x_m^1,\ldots ,x_m^k\right) \right\| \leq
\lim_m \sup_n \left\| A\left( x_n^1,\ldots ,x_n^{k-1},x_m^k\right)\right\|
= 0\, .
$$
In the general case, the proof follows that of Theorem~6 in
\cite{GGwh}.

(b) $\Ra$ (c) is obvious.

(c) $\Ra$ (a) is clear.
\finesp

The condition ${\cal L}(E,c_0)=\CC (E,c_0)$ implies
that $E$ has the Dunford-Pettis property.
However, there are spaces with the Dunford-Pettis
property that admit noncompletely
continuous operators into $c_0$ (e.g.\ $E=c_0$, $E= L_1[0,1]$).

\vspace{\fill}\hspace*{\fill}{\small file pgrp.tex}

\end{document}